\documentclass[reqno,11pt]{amsart}
\usepackage{amscd,amssymb, amsmath}
\usepackage[pdftex]{graphicx}
\usepackage[cp1251]{inputenc}
\usepackage[T2A]{fontenc}
\usepackage[russian, english]{babel}
\usepackage{mathabx}

\def\acts{\lefttorightarrow}

\newtheorem{theorem}{Theorem}[section]
\newtheorem{lemma}[theorem]{Lemma}

\theoremstyle{definition}
\newtheorem{definition}[theorem]{Definition}

\newtheorem{claim}[theorem]{Claim}

\newtheorem{question}[theorem]{Question}

\newtheorem{note}[theorem]{Note}

\begin{document}
\title{An operad from the secondary polytope}
\author{Lev Soukhanov}
\maketitle

\begin{abstract}
We study the enumerative geometry of orbits of multidimensional toric action on projective algebraic varieties and develop a new cyclic differential-graded operad, conjecturally governing the real version of the enumerative geometry of these toric orbits.
\end{abstract}

\section{Introduction}
The description of the space of orbits of multidimensional toric action $(\mathbb{C}^*)^d \acts \mathbb{P}^n$ was a part of a hypergeometric project in [GKZ] and led to a development of a notion of a \textbf{secondary polytope}. The space of orbits (together with the limit degenerations) is a toric variety which is dual to the \textbf{secondary polytope}. In [GKZ] this space is called the \textbf{Chow quotient}.

In [A], this space was broadly generalised to the moduli space consisting of multiple irreducible components, called the \textbf{modular Chow quotient}. Virtual fundamental classes for this space are defined in the work [M].

In our work we construct from the combinatorics of secondary polytope the \textbf{Secondary operad}: cyclic operad with operations corresponding to the strata of the real version of modular Chow quotient for $d=2$. It corresponds to the Morse theory in the $1$-dimensional setting, and the complex version (which is still unclear) would correspond to the Losev-Manin commutativity equations [LM],[L].

The relations between the various problems of enumerative geometry discussed in this paper are presented in the following tables:

\begin{center}
\begin{tabular}{|c|c|c|c|}
\hline
Theory: & Morse & $\mathbb{C}^*$ action & enumeration\\
 & & & of curves\\

\hline
Enumerative problem: & integral & invariant & Gromov-Witten\\
  & curves & rational curves & \\
\hline
Moduli space: & A-$\infty$ operad & Losev-Manin & Deligne-Mumford\\
\hline
Algebraic structure: & $\partial^2 = 0$ & commutativity & WDVV\\
 & & equations & \\

\hline
\end{tabular}
\end{center}

The relation between the commutativity equations and the Morse theory is explained in [L], and relation between commutativity equations and WDVV is explained in [LM]. We'd like (ideally) to construct the analogue of the Gromov-Witten theory for the mappings of the surfaces to an algebraic variety $X$. There are various difficulties, but one of the ways of doing it might be looking for its equivariant analogies.

\begin{center}
\begin{tabular}{|c|c|c|c|}
\hline
Theory: & 2-Morse* & $(\mathbb{C}^*)^2$ action & enumeration\\
& & &of surfaces\\
\hline
Enumerative problem: & integral leaves & ? & ?\\
\hline
Moduli space: & ? & Chow quotient & Alexeev's\\
& & & stable pairs\\
\hline
Algebraic structure: & Secondary operad & ? & ?\\
\hline
\end{tabular}
\end{center}

* - By $2$-Morse we mean the description of integral leaves of the pair of commuting gradient-like vector fields. This theory seems to lack argument of general position (which is present in $1$-Morse), so it seems one should think only of the complexifiable actions. Nevertheless, we can still develop the corresponding algebraic structures.

The "?" signs in the table above correspond to the fact that we don't really understand, what kind of enumerative problems should be considered: should we enumerate surfaces passing through the number of points, or the number of curves. The construction of Alexeev suggests the latter, but precise problem is unclear.

However, "Morse-like" real version of the problem can be considered without any deformation, and this advantage allows us to describe the structure for the real case.

The paper is organised as follows. In Section $2$ we will construct the formula for the expected codimension of a toric degeneration, show that it sometimes gives non-positive result and suggest a program of getting rid of such "bad" degenerations. Section $3$ is a brief exposition of the theory of secondary polytopes and Alexeev's modular Chow quotient. In the Section $4$ we will define perturbed regularity condition and use it to fulfill our program and define the secondary operad. In the last section we will discuss questions our work opens.

The algebro-geometric part of the paper is mostly sketchy, because we don't prove anything new in the algebro-geometric realm. The reader interested in the rigorous constructions of the moduli spaces of surfaces should definitely take a look in the works [A],[M]. However, we hope that the algebraic structures from this enumerative geometry after development would incorporate these formalisms and show further directions for the appropriate generalisations.

\textbf{Acknowledgements.} I want to thank my scientific advisor, Andrei Losev, for presenting me this question and constantly keeping my interest to the subject, even when it seemed to be completely obscure for both of us. I want to thank Petr Pushkar and Dmitry Korb for the discussions and counterexamples strongly influencing my intuition in this area and Valery Alexeev, who's minicourse on the moduli spaces of surfaces at HSE gave me the more rigorous basis for thinking about the geometric counterpart of the problem. Also I want to thank Anastasiya Goryatchkina, who helped me with vector graphics for this paper.

  The article was prepared within the framework of a subsidy granted to the HSE by the Government of the Russian Federation for the implementation of the Global Competitiveness Program. The author was also supported in part by the Moebius Contest Foundation for Young Scientists and by the Simons Foundation.

\section{Expected codimension}

Let us consider a real smooth compact manifold $M$ with a gradient-like vector field $v$ in a general position. Let us consider two critical points $A$, $B$. $n_A^B$ is defined as a number of integral curves of the field $v$, starting at $A$ and ending in $B$, taken with appropriate signs (if their number is infinite, $n$ is taken to be $0$). The main theorem of the Morse theory states that $n^2 = \sum \limits_{B} n_A^B n_B^C = 0$.

Let us present the sketch of proof (or, maybe, explanation) of this theorem. For the detailed exposition of the Morse theory the reader can take a look in [Mil].

Let us consider the space of integral curves between $A$ and $C$. If it is $k$-dimensional, then its boundary is $k-1$-dimensional and consists of degenerated curves. The degeneration in the general position (i.e. degenerations which occur in the codimension one) are the degenerations with one intermediate critical point ($B$). Taking $k=1$ we can conclude that the space of curves between $A$ and $C$ consists of the union of circles and line segments, and boundary components of the line segments are precisely the elements of the sum $\sum \limits_{B} n_A^B n_B^C$. For each line segment, the contributions of its ends have an opposite signs, hence the sum vanishes.

If we consider the fundamental cycle of the space of integral curves $\Phi_A^C$, then by the same argument, the following Maurer-Cartan relation holds: $\partial (\Phi_A^C) = \sum \limits_{B} \Phi_A^B \Phi_B^C$ for any $A$, $C$.
\begin{center}
\includegraphics[width=6.0cm, viewport=100 280 550 730, clip]{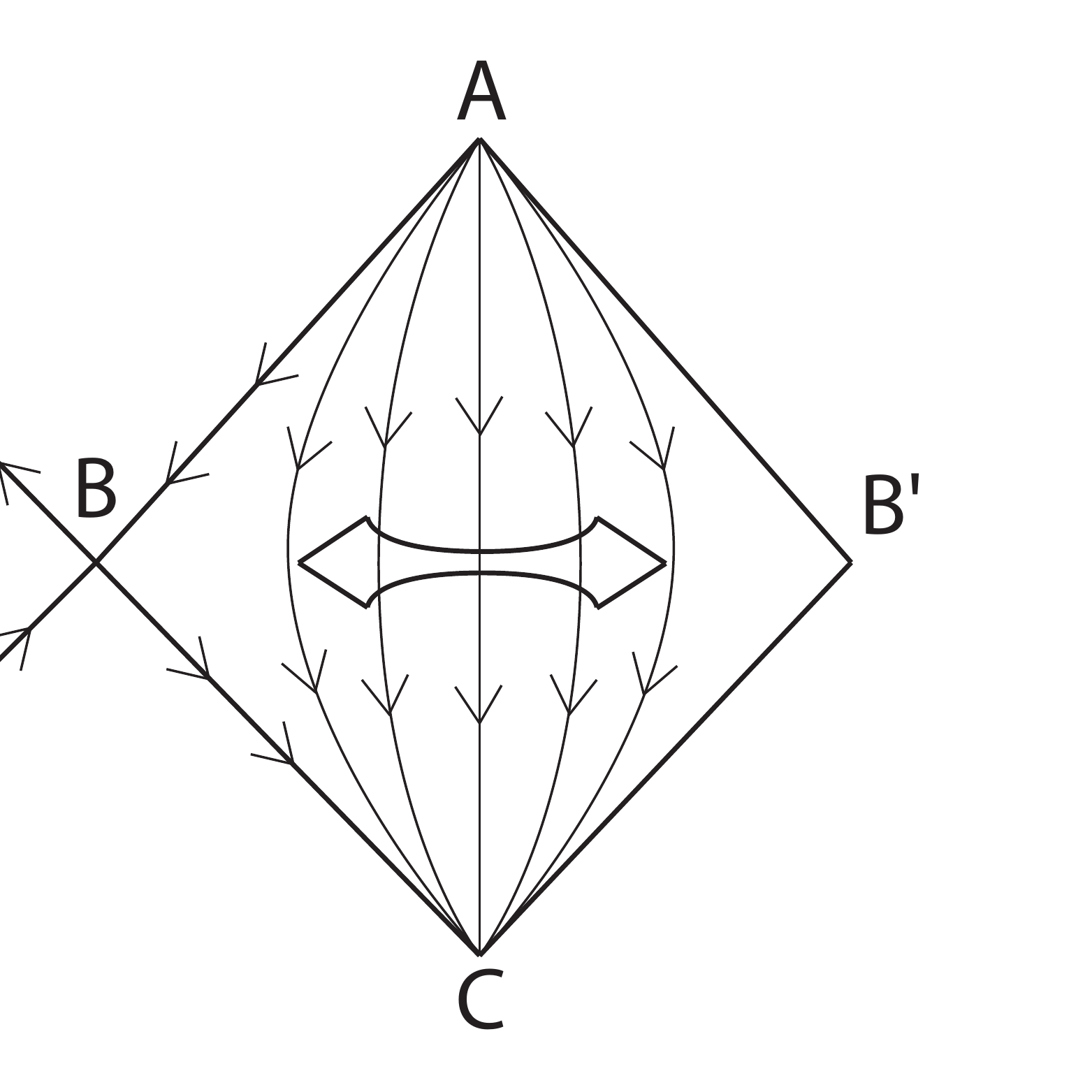}

\textit{Degeneration of an integral curve in a Morse theory.}
\end{center}
The same program can (and should) be applied to the leaves of an action of a pair of commuting vector fields (or, in the complex case, to the orbits of $(\mathbb{C}^*)^2$-action).

The closure of an orbit of $(\mathbb{C}^*)^d$ is a toric variety (and in the real case, closure of the leaf is a real toric variety). So, at first we will describe an applicable combinatorial formalism for the degenerations.

\definition \textbf{Affine fan} $F$ is the finite decomposition of the space $\mathbb{R}^d$ into the union of convex (maybe not bounded) polytopes. The fan is always assumed to be non-empty.

Affine fan encodes the following variety: let us take a vertex $v$ of a fan $F$. Link of this vertex is a toric fan $T(v)$. Let us take the union of $T(v)$'s and glue* them along the components of the boundary, corresponding to the edges of the fan $F$.

* - gluing in algebraic geometry can be done in multiple ways, but there is the most "economic" one, which is seminormal.

\definition \textbf{Infinity of the affine fan} $Inf(F)$ is a proper toric fan, constructed from the affine fan as follows: take any point in the space $\mathbb{R}^d$ and start rescaling it with constant tending to zero. Limit of this procedure is called infinity of the affine fan.

Affine fan encodes the flat equivariant degeneration of a toric variety, encoded by its infinity.

\begin{center}
\includegraphics[width=6.0cm, viewport=150 250 500 600, clip]{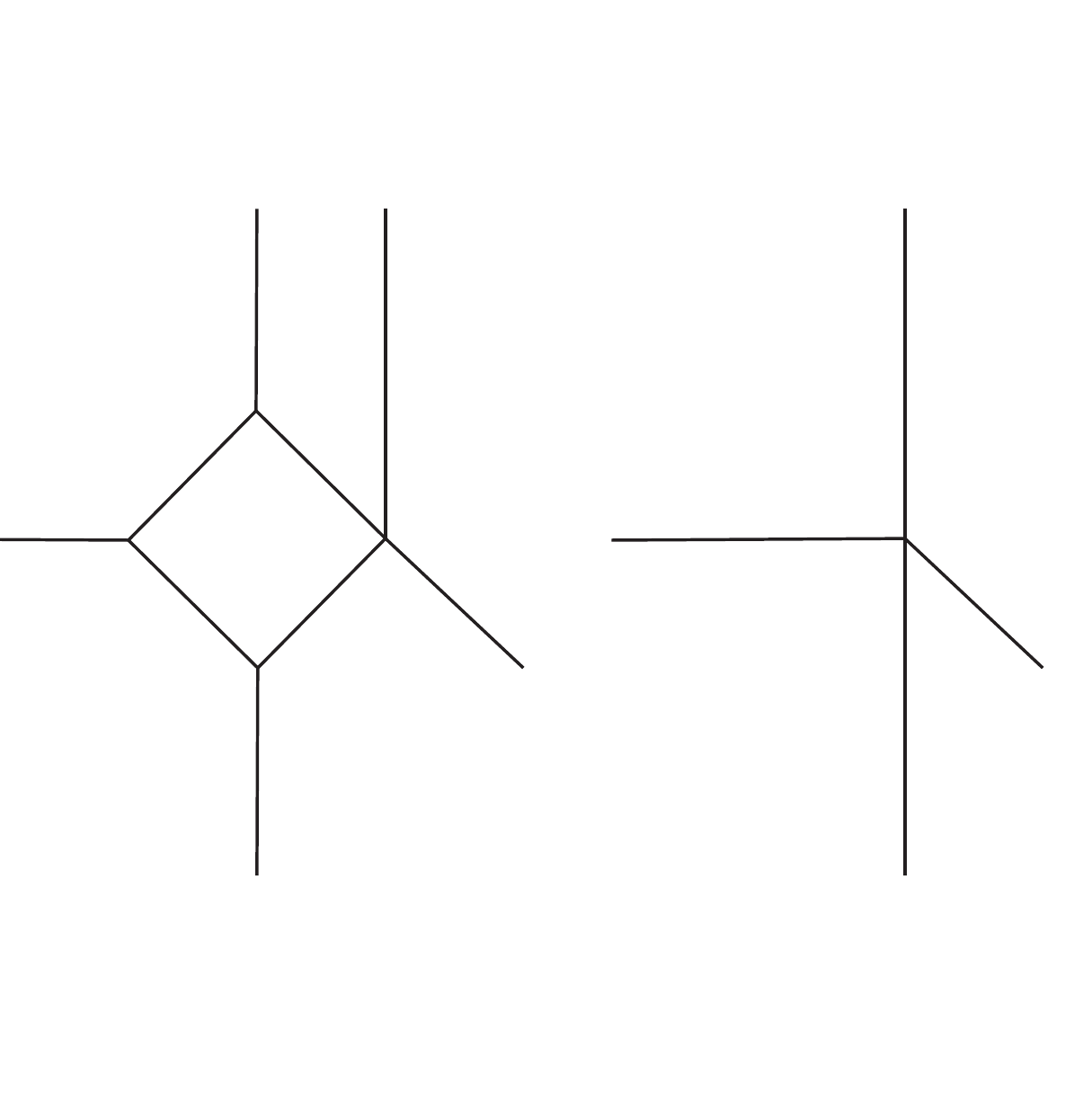}

\textit{Affine fan and its infinity.}
\end{center}

Let us consider the $2$-dimensional affine fan $F$. Let us denote by $v(F)$ the number of its vertices, and $e(F)$ the number of its internal edges.

Let $X$ be an algebraic variety with $(\mathbb{C}^*)^2$-action (corr. its real part).

\definition \textbf{Expected codimension} of $F$ is a jump of an invariant euler characteristic of a normal sheaf to the general leaf and the degeneration of type $F$.

\lemma It depends only on $F$ and calculated as

$Codim(F) = 2 v(F) - e(F) - 2$.

\proof

The restriction of a tangent bundle of $X$ is flat, so the jump of an invariant euler characteristic of a normal sheaf can be calculated as a minus jump of an invariant euler characteristic of a tangent bundle. It proves the first statement.

An invariant euler characteristic of the tangent bundle is an expected dimension of equivariant automorphisms - which is $2 \times$ number of toric components $-$ number of relations, i.e. internal edges. The non-degenerated leaf has a $2$-dimensional space of automorphisms (the torus itself), hence the jump equals $2 v(F) - e(F) - 2$. $\blacksquare$

\begin{center}
\includegraphics[width=12.0cm, viewport=-200 400 800 800, clip]{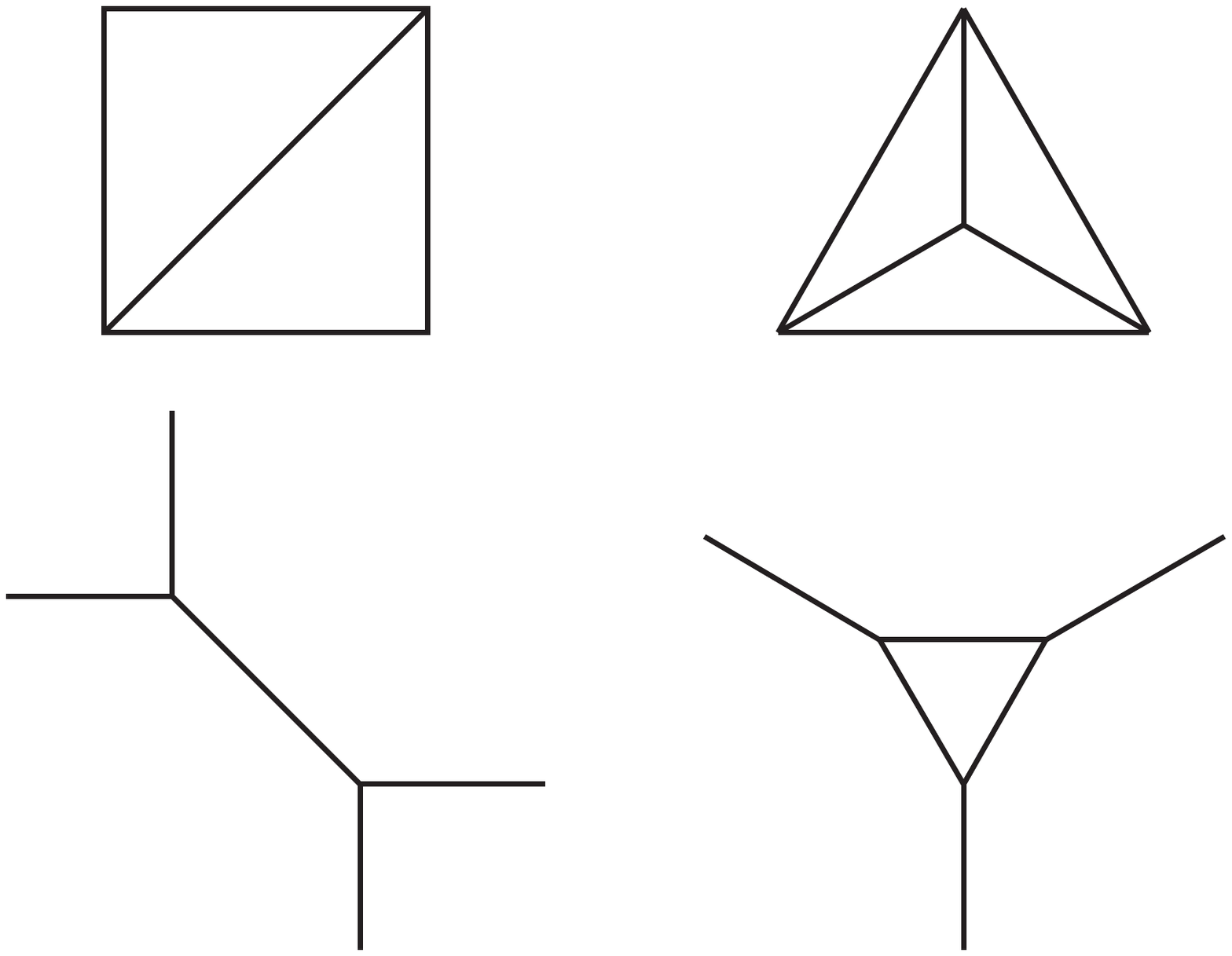}

\textit{Various degenerations of codimension 1. Diagrams in the upper part of the picture show the combinatorics of the leaves.}
\end{center}

So, by the analogy with the Morse theory, one could possibly assume the following \textbf{naive conjecture}:

Let us consider $M$ - real smooth compact manifold with a \textit{good} action of two commuting gradient-like vector fields. Let us consider, additionally, that $0$-dimensional and $1$-dimensional orbits of the group $\mathbb{R}^2$ of flows of these vector fields are isolated. Then, for any cycle $C$ of $1$-dimensional orbits we can consider the number of integral leaves $n_C$.

\textbf{Naive conjecture}: Then, the following equality holds (the sum is taken with appropriate signs) for any cycle $C$ of $1$-dimensional orbits.

\begin{equation}
\sum \limits_{F|codim(F)=1} \#F(C)
\end{equation}

where $\#F(C)$ is a number of degenerated toric varieties of type $F$ with the boundary $C$.

This naive conjecture is not true and cannot be true in any reasonable way and the simplest way to see it is the following:

\textit{Expected codimension sometimes happen to be zero or even negative}.

The example is given below.

\begin{center}
\includegraphics[width=10.0cm, viewport=0 300 550 730, clip]{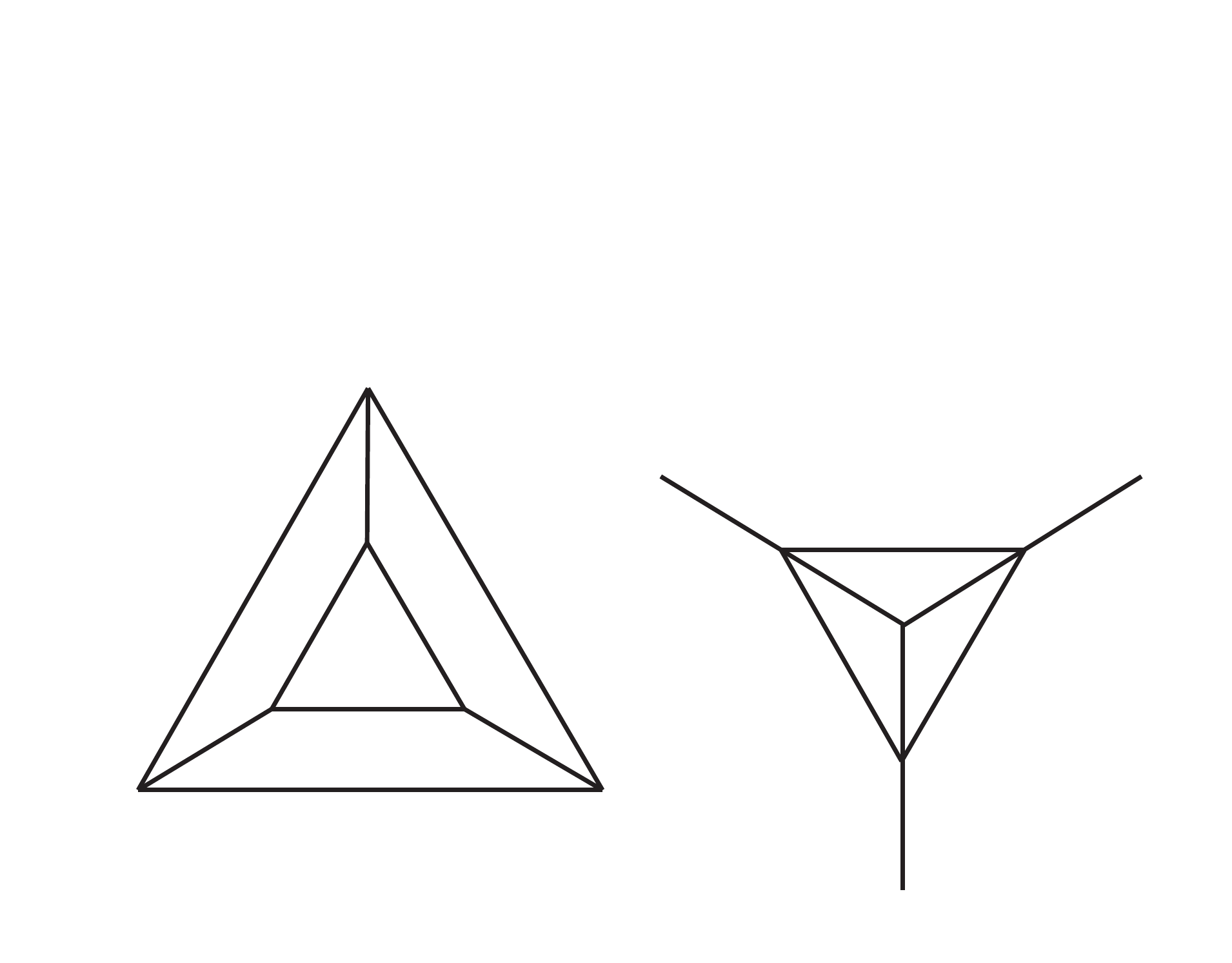}

\textit{Zero codimension degeneration.}
\end{center}

In the next section we will see that this situation is typical and the degeneration \textit{actually has} more deformations than the non-degenerated leaf. These deformations belong to the other irreducible component of the moduli space, constructed by Alexeev.

Let us denote by $\Phi_C$ the fundamental class of the space of leaves with the boundary $C$ and by the $F(\Phi_{C_1}, ..., \Phi_{C_n})$ the fundamental class of the space of degenerated leaves of type $F$ with components with boundaries $C_1, ..., C_n$.

\textbf{Idea of the treatment:} An equation $\partial(\Phi_C) = \sum \limits_{F|codim(F)=1} F(\Phi_{C_1}, ..., \Phi_{C_n})$ cannot hold, because there are also boundary components, corresponding to the $codim(F) < 1$. But if we consider, for example, $F(\Phi_{C_1}, ..., \Phi_{C_n})$ for $codim(F) < 1$, we will get the class of wrong dimension. We need only the class of the locus of such degenerated toric varieties, which occur as the degenerations of the regular leaves. The idea is to deform this locus to the homotopy-equivalent component of the boundary.

\section{Secondary polytopes and Alexeev's modular Chow quotient}

This part of exposition mainly follows [GKZ].

Let us consider a set $A$ of points of a lattice $\mathbb{Z}^2$. Denote by $X$ the projective space $\mathbb{P}(\mathbb{C}^A)$ together with the action of $(\mathbb{C}^*)^2$, defined as follows: for a basis vector of $\mathbb{C}^A$, corresponding to the point $a \in A$ with coordinates $(x, y)$ it acts with the character $z^x w^y$.

The closure of the general orbit $Z \in X$ is a $2$-dimensional toric variety, which moment polytope coincides with the convex hull of $A$.

\definition Irreducible component of $Z$ in the $(\mathbb{C}^*)^2$-invariant locus of the Chow scheme is called \textbf{Chow quotient} of X.

\note It is a toric variety w.r.t to the torus $(\mathbb{C}^*)^A/(C^*\times(\mathbb{C}^*)^2)$ which naturally parametrises general orbits and their possible degenerations.

\definition Connected component of $Z$ in the $(\mathbb{C}^*)^2$-invariant locus of the Chow scheme is called \textbf{modular Chow quotient} of X.

\note This definition is different from the original definition in [A], as it doesn't take into account the stack structure of the modular Chow quotient (which we won't use).

Now we shall construct a toric fan of the Chow quotient.

Let us consider a space $Fun(A) = \mathbb{R}^A$ of real-valued functions on $A$.

The space of affine functions $Aff(\mathbb{R}^2)$ is naturally mapped to $Fun(A)$ via restriction. We denote by $V = Fun(A)/Aff(\mathbb{R}^2)$.

For any $v \in Fun(A)$ let us consider a minimal concave function $f$ on the convex hull of $A$ such that $f(x,y) \geq v(a)$ for $\forall a \in A$ with coordinates $(x, y)$. It is obviously piecewise-linear, with the domains of linearity being convex polytopes $\{P_1, ..., P_k\}$, $\bigcup P_i = Conv(A)$.

It depends only on a class of $v$ in $V$.

Denote the corresponding decomposition of $Conv(A)$ into the union of polytopes by $D(v) = \{P_1, ..., P_k\}$

This data is combinatorial and induces the stratification of a space $V$ into a union of open convex polyhedral cones.

\definition This collection of cones is called the \textbf{secondary fan} of $A$.

\note Elements of $V$ naturally encode $1$-parametric orbits in the torus $(\mathbb{C}^*)^A/(\mathbb{C}^*\times(\mathbb{C}^*)^2)$, and $D(v)$ encodes the corresponding degeneration of a (moment polytope of a) leaf while moving along this $1$-parametric orbit.

\definition The decomposition $D = \{P_1, ..., P_k\}$ is called \textbf{regular} if it comes from the construction as $D(v)$ for some $v$, i.e. if it corresponds to some degeneration of a general orbit.

\note In the modular Chow quotient, the leafs of the degeneration deform independently, hence no condition of regularity is imposed.

\definition An affine fan $F$ in the space $(\mathbb{R}^2)^*$ is called \textbf{normal} to the decomposition $D$ in the space $\mathbb{R}^2$ if its vertices correspond to the polygons of $D$ and the edges are normal to the boundaries between corresponding polygons.

\begin{center}
\includegraphics[width=6.0cm, viewport=0 200 550 500, clip]{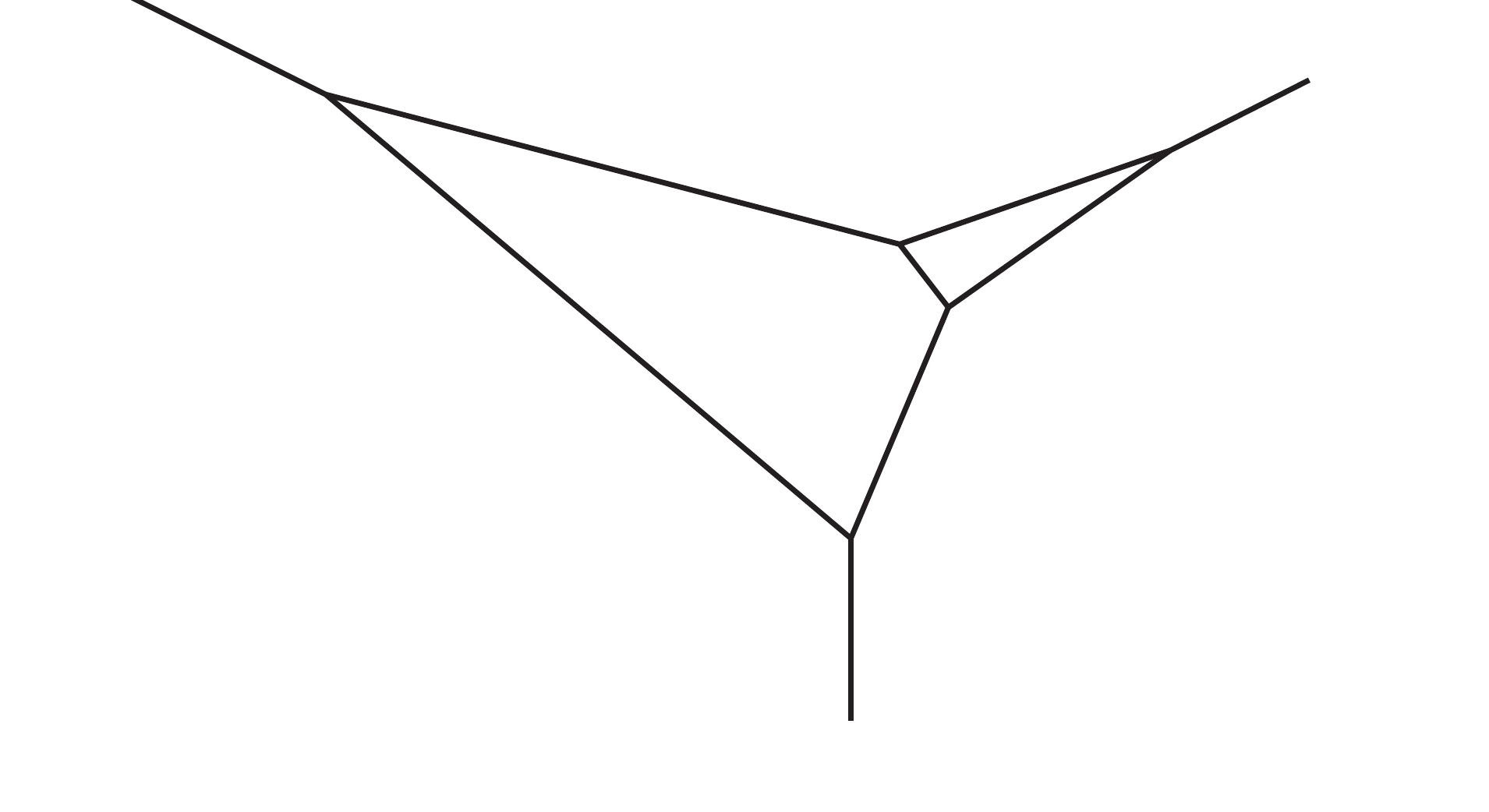}
\includegraphics[width=6.0cm, viewport=100 230 550 730, clip]{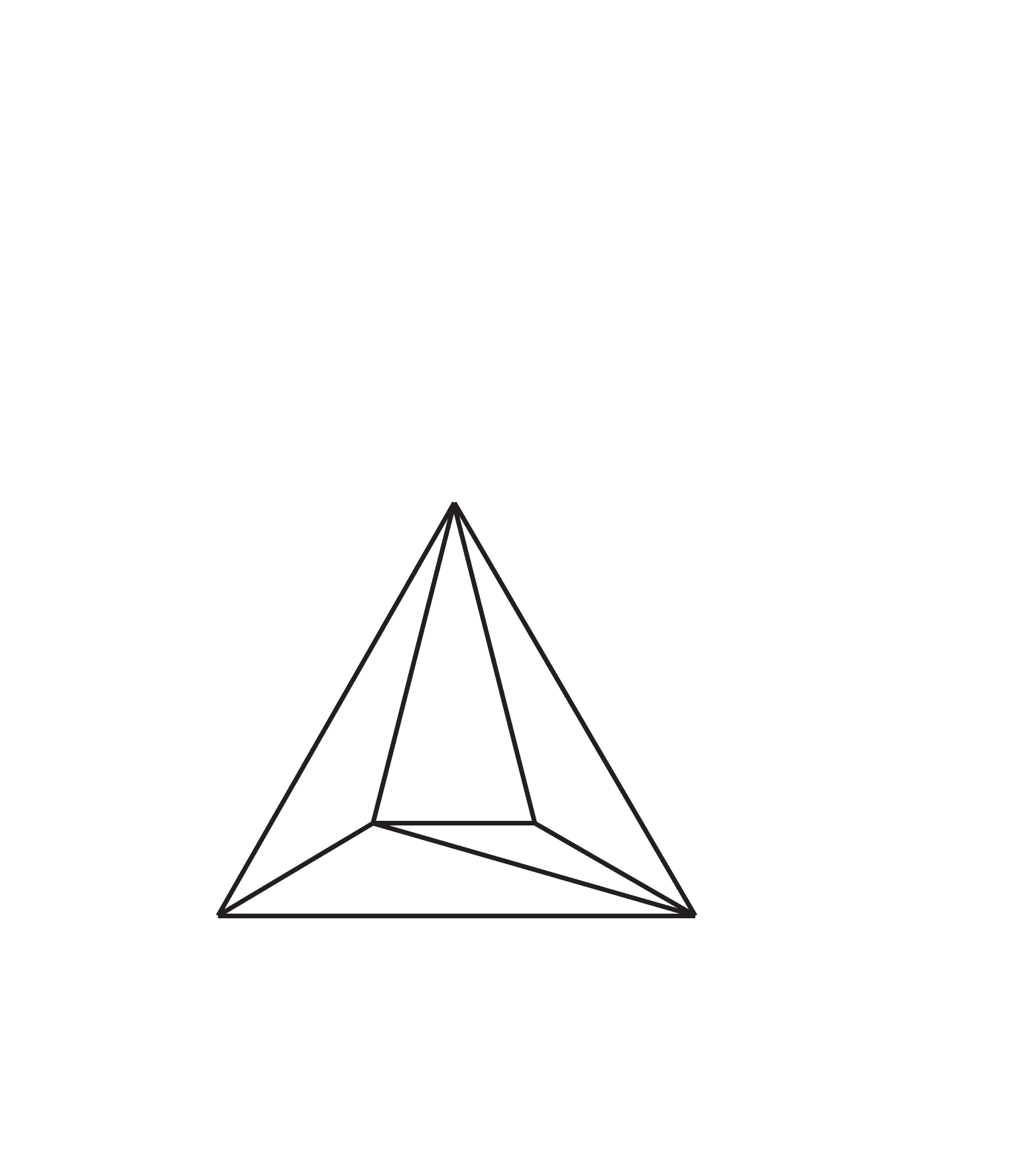}

\textit{Normal fan to the decomposition.}
\end{center}

\theorem Normal fan to the decomposition $D$ exists if and only if this decomposition is regular.

\proof

Let us assume that the decomposition is regular w.r.t. to the concave function $f$. Consider affine functions $f_i = f|_P$ on each $P_i$. Turn them into linear functions $\gamma_i \in (\mathbb{R}^2)^*$. It will be the set of vertices of the affine fan. For two neighboring polygons $P_i$ and $P_j$ consider an edge between $\gamma_i$ and $\gamma_j$. It will be orthogonal to their common boundary, because $f_i = f_j$ on it. It is not hard to see that concavity of function $f$ guarantees that the resulting object will be an affine fan.

Conversely, having an affine fan, we get the set of $\gamma_i$'s for free. Let $g_i$'s be any affine functions with the linear parts $\gamma_i$'s. Then, $g_i - g_j|_{P_i\cap P_j}$ are constants $a_{ij}$, satisfying the cocycle condition $a_{i_1 i_2} + ... + a_{i_s i_1} = 0$ for any $P_{i_1}, ..., P_{i_s}$ located along a vertex in a decomposition. $H^1(Conv(A), \mathbb{R}) = 0$, hence this cocycle is a coboundary, so we can produce another set $f_i = g_i + const_i$ such that $f_i = f_j$ on $P_i \cap P_j$. Concavity of the constructed piecewise-linear function follows (by the local checks on a boundaries between $P_i$ and $P_j$'s) from the fact that $\gamma_i$'s are the vertices of an affine fan. $\blacksquare$

So, we are in a need for the new object, which would be dual to the non-regular decompositions. From the fact that two affine fans which are combinatorially equivalent and have parallel edges correspond to the same degenerated toric variety we come to the following

\definition \textbf{Graph with directions (GwD)} is an oriented graph with outgoing edges (one can think that they go to the vertex "infinity"). Each edge is endowed with an additional data: the direction, i.e. unit vector in $\mathbb{R}^2$. It is allowed to change orientation of an edge together with changing the direction of an edge to the opposite.

We will use the wording "outgoing direction along an edge $e$ from a vertex $v$" for the direction of $e$ if it is oriented from $v$ and for the minus direction of $e$ otherwise.

We will consider only graphs such that outgoing directions from any vertex positively generate $\mathbb{R}^2$.

\note It is the first (but not the last) place in the paper when giving a definition in the dimension $d > 2$ would be hard (and we don't even know how to do it right).

\definition \textbf{Positive representation of a GwD} $\Gamma$ is an affine fan $F$ with vertices corresponding to the vertices of $\Gamma$, edges corresponding to the edges of $\Gamma$ and parallel to the directions of those edges (orientation-preserving).

They are assumed to be equivalent if they differ only by a parallel transport.

\begin{center}
\includegraphics[width=10.0cm, viewport=-100 300 700 600, clip]{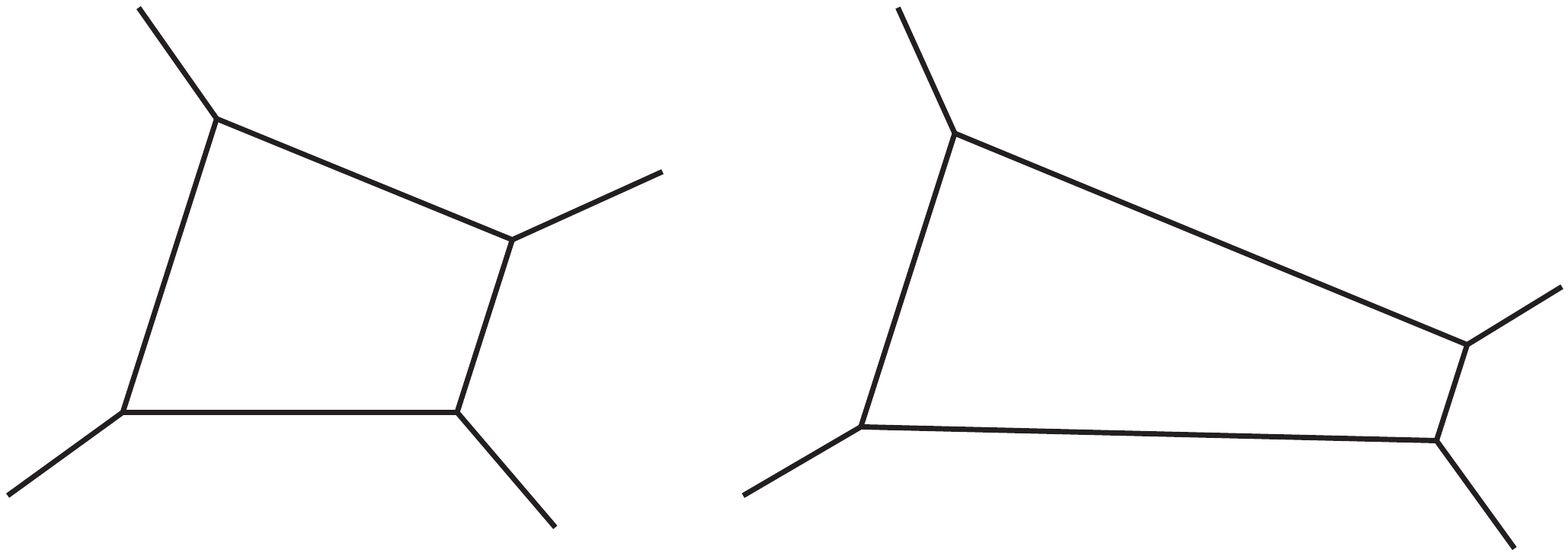}

\textit{Different positive representations of the same GwD.}
\end{center}

\definition \textbf{Representation of a GwD} $\Gamma$ is a set of points $p_i$'s corresponding to the vertices of $\Gamma$, such that vectors $e_{ij} = p_i - p_j$ are parallel to the directions of edges between vertices $i$ and $j$. No conditions about orientation are imposed, $p_i$'s can even coincide (and edge between such points would be parallel to anything).

\note Representations of $\Gamma$ obviously form a vector space $Rep(\Gamma)$, and set of positive representations $PRep(\Gamma)$ is a convex polyhedral cone.

The reader could notice that we actually defined some invariant of the GwD - the dimension of its space of representations. The calculation of dimensions shows that it is expected to be $2v - e - 2$, where $v$ denotes the number of vertices and $e$ denotes the number of internal edges - it coincides with the formula of the expected codimension.

\definition GwD $\Gamma$ is called \textbf{too rigid} if $dim(Rep(\Gamma))$ is greater than expected.

\note Degenerations corresponding to the too rigid GwD's are the boundary components, connected to the other irreducible components of modular Chow quotient. It is (maybe in not so explicit way) can be found in [A].

\note Those GwD's wouldn't exist for a general set of directions. However, we cannot simply deform the directions - because they are rational edges of toric fans.

\begin{center}
\includegraphics[width=8.0cm, viewport=0 350 550 700, clip]{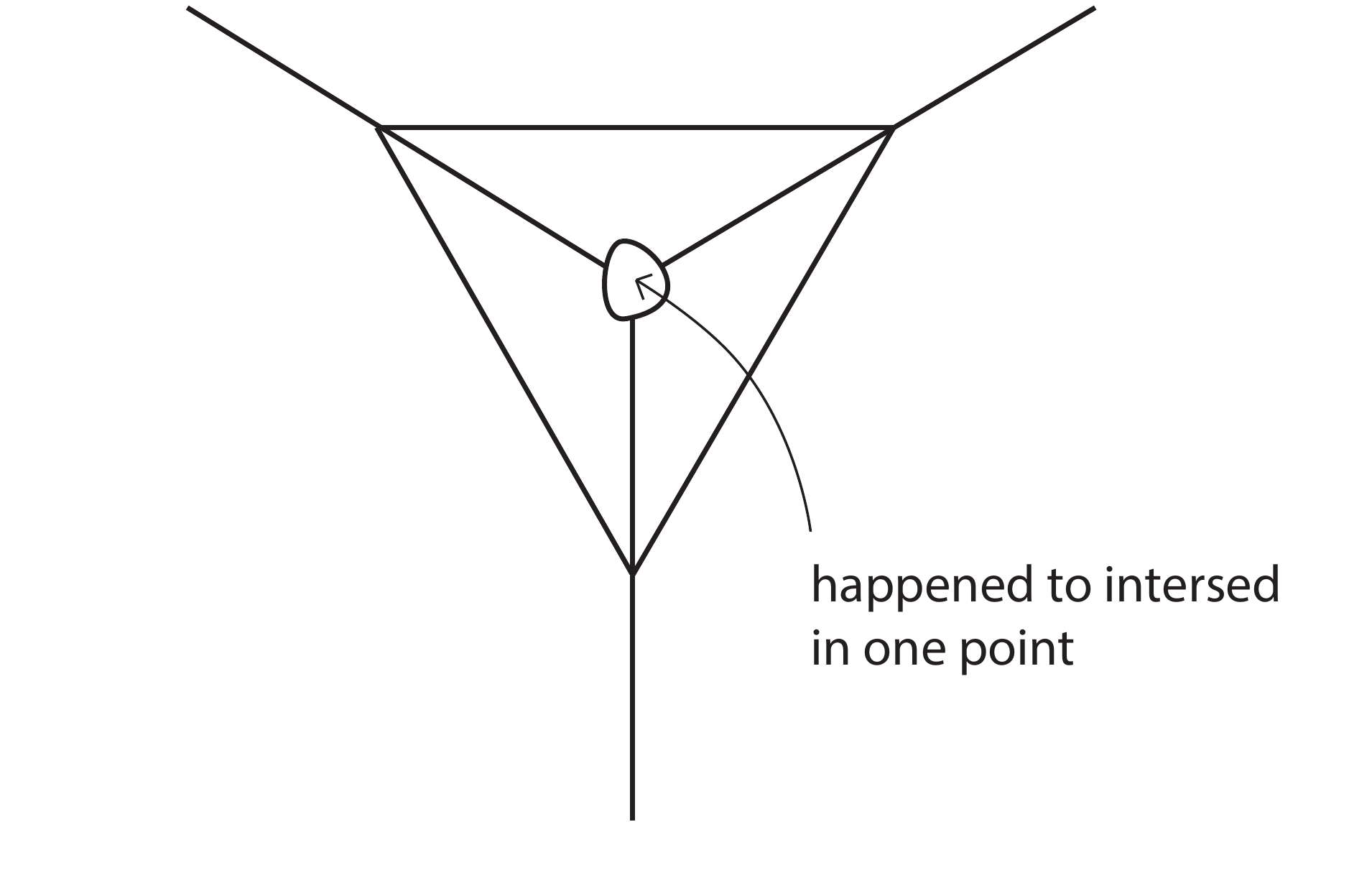}

\textit{Too rigid graphs doesn't exist for a general set of directions.}
\end{center}

We proceed by construction of the secondary operad.

\section{Secondary operad}

Let $A$ be a (finite) set in $\mathbb{R}^2$. Assume that there are no three points in $A$ on the same line (this condition is technical and will be discussed later).

\definition The basis operations of a \textbf{secondary operad} are any polygons with vertices in $A$, decomposed into a union of convex polygons with vertices in $A$ with non-intersecting interiors. The composition is the gluing of two polygons by a set of edges.

\textbf{Sign convention:} the set of internal edges of the polygonal decompositions is ordered. Changing an order shifts the sign by the sign of the permutation. The parity of the basis element is counted as the parity of the number of edges plus one, hence the composition is of degree $0$. There might be different conventions about $\mathbb{Z}$-grading.

\note For a regular decomposition $D$ the sign data is equivalent to the orientation on the polyhedral cone $PRep(F)$, where $F$ is a normal fan to $D$ (coming from the exact sequence $0 \rightarrow Rep(F) \rightarrow (\mathbb{R}^2)^v \rightarrow (\mathbb{R})^e \rightarrow 0$). Hence, any regular decomposition of codimension $1$ has a canonical sign data, as it has the $1$-dimensional $PRep$.

\textit{Assume that all GwDs dual to the decompositions of} $A$ \textit{are not too rigid (i.e., all decompositions are regular).}

\note It is rare situation and it doesn't follow from the conditions we imposed on $A$.

Then, the following works:

\definition The differential $\partial$ is defined as follows: for any convex polytope $P$ it returns $\partial(P) = ( \sum \limits_{D|_{codim(D) = 1}}D)$ the sum over decompositions of $P$ (with the canonical sign data). For a decomposition $\partial(D)$ is defined by Leibniz rule.

\lemma $\partial^2 = 0$

\proof The sign data from the Leibniz rule is equivalent to the orientation by outward normal vector: for any decomposition $D$ of codimension 2 there is $2$-dimensional cone of representations, which has boundary, consisting of two edges - corresponding to the two summands of type $D$ in $\partial^2 (P)$ with the opposite signs.

\note Reformulating geometrically: consider the space of representations of GwD's dual to the decompositions of $P$, and glue them along the boundary components - when a part of the representation of a GwD shrink into a point, glue this boundary component to another space of representations with GwD with this part shrinked. Then the operations of secondary operad will correspond to the (oriented) facets of this space, $\partial$ is just the dual cellular differential and $\partial^2 = 0$ is obvious.

\note The space constructed above is actually a toric fan for a real version of Chow quotient.

\begin{center}
\includegraphics[width=12.0cm, viewport=0 400 650 750, clip]{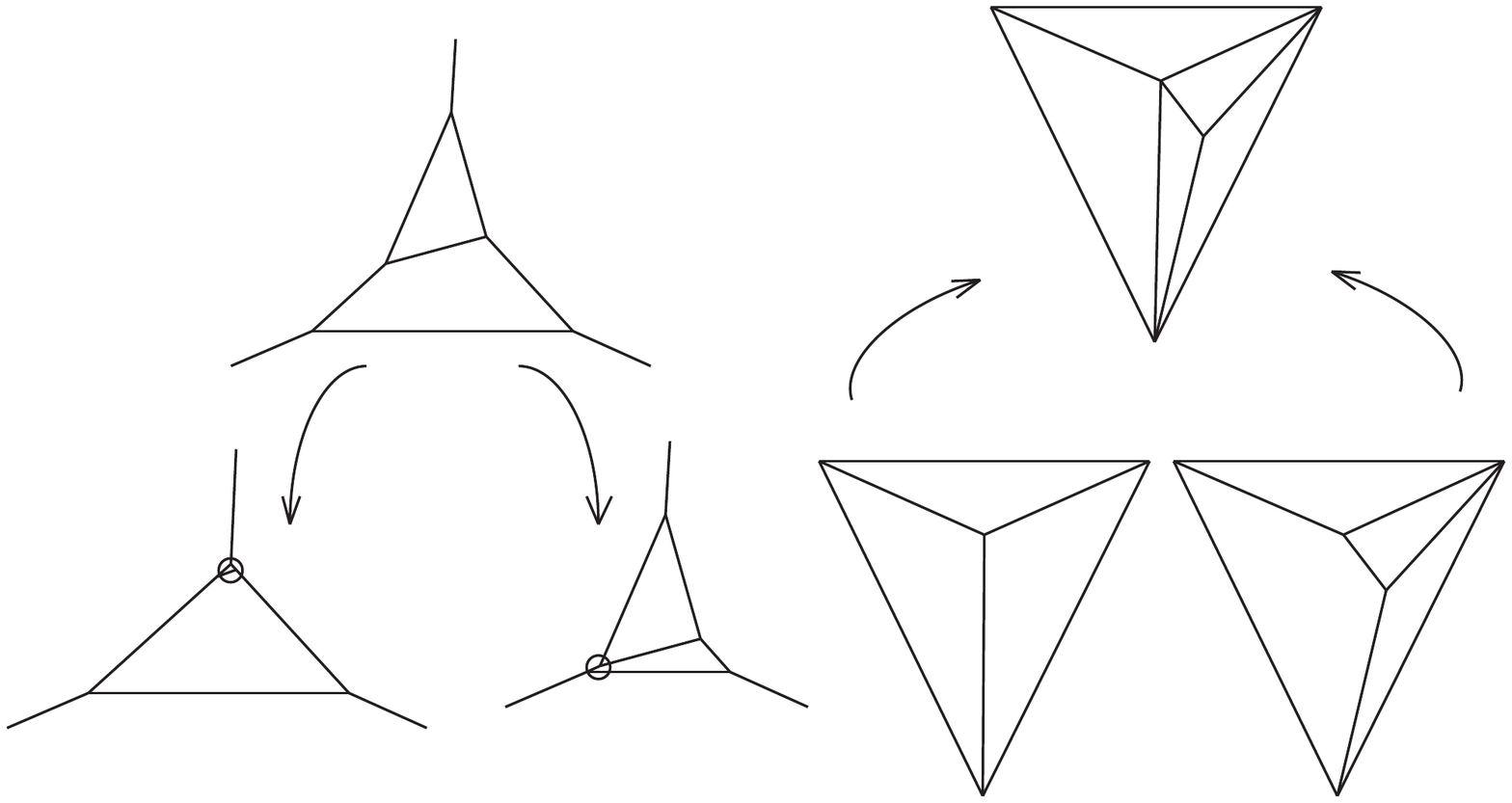}

\textit{Shrinks of the codim=2 GwD and $\partial^2 = 0$.}
\end{center}

Now we finally move to the case where too rigid graphs occur.

Let us consider the set of all possible directions $S$ occuring in all decompositions of $A$. Let us consider the (germ of) deformation $S_t$ general enough (each direction is being deformed smoothly). Then, no too rigid graphs can be constructed using this set of directions for $t \neq 0$.

\note This is the main reason why we ask $A$ to have no triples of points on a same line - we don't want the same directions to occur in a different parts of GwD.

\definition The decomposition $D$ is called \textbf{perturbedly regular} if the dual GwD with deformed directions admits a positive realisation for $t \neq = 0$.

\note No too rigid graphs are dual to the perturbedly regular decompositions, but some irregular decompositions will become perturbedly regular.

\begin{center}
\includegraphics[width=16.0cm, viewport=0 600 700 800, clip]{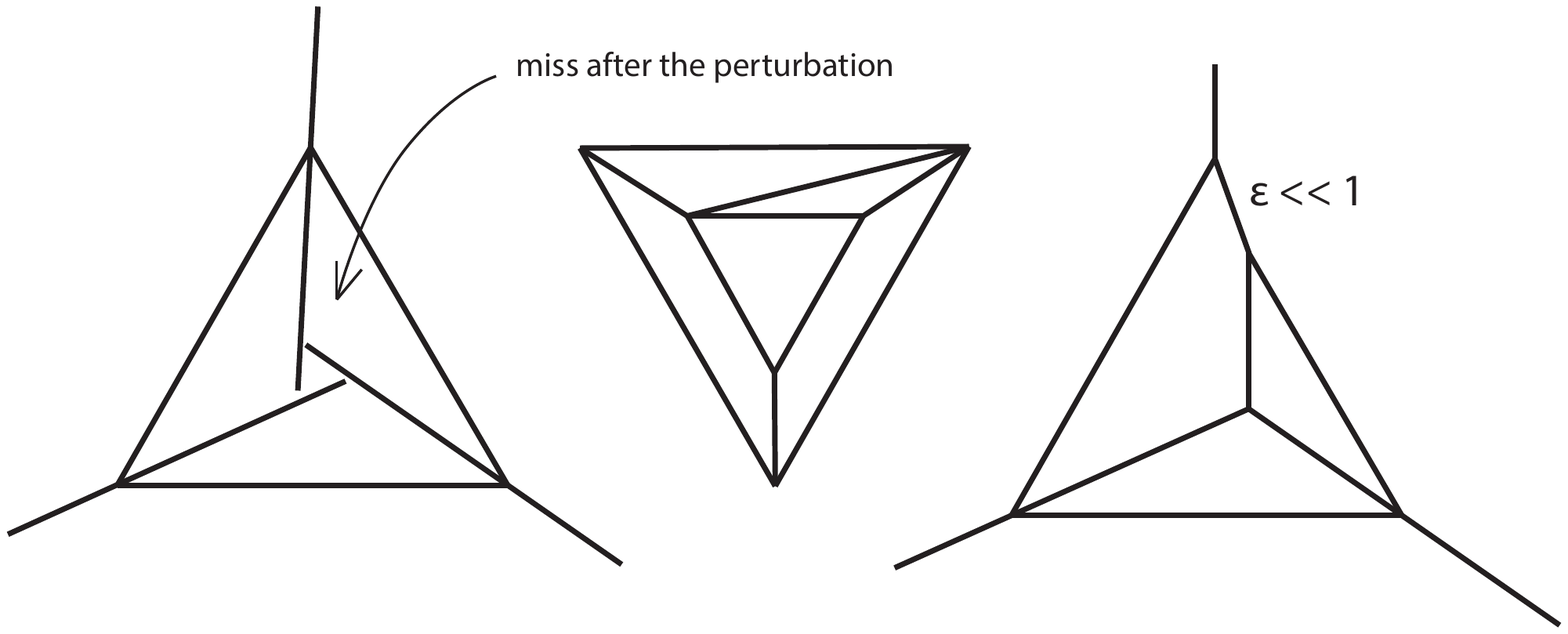}

\textit{Some irregular decompositions will become perturbedly regular.}
\end{center}

\definition $\partial(P) = ( \sum \limits_{D|_{codim(D) = 1}}D)$, where the sum is taken over all perturbedly regular decompositions of codimension 1.

\theorem $\partial^2 = 0$

\proof The same argument as in the previous theorem should be applied to the same space with $0 < t \ll 1$.

Now, as we defined an operad, we will discuss why do we think that this definition is right.

Consider the following situation: the leaf deforms to the degeneration corresponding to the too rigid graph. Then these degenerations form a locus in an irreducible component of the modular Chow Quotient. The fundamental class of this locus should be the summand in the $\partial(\Phi)$, but we want it to be the composition in an operad, so we deform it to the union of boundary components of this component.

\definition \textbf{Perturbation} of a too rigid graph is a set of all perturbed graphs who's realisations tend to this too rigid graph as $t \rightarrow 0$.

\note So, the definition of $\partial$ can be reformulated as follows: we take the sum over all $F$'s with $codim \leq 1$, but then perturb too rigid summands.

\claim {Let us take the real part of the modular Chow quotient. Consider too rigid decomposition $D$ and the corresponding irreducible component $Z$ (which is just the product of modular quotients for each $P_i \in D$). Then the intersection $L$ with the main component can be deformed (without moving a boundary) to the union of the boundary components from the perturbation of $D$.}

\textit{Sketch of the proof}

$Z$ is a toric variety w.r.t. to the torus $\prod \limits_i(\mathbb{C}^*)^{P_i}/((\mathbb{C}^*)^2\times(\mathbb{C}^*))$, and $L$ is a closure of an orbit of the subtorus which is given by the image of $(\mathbb{C}^*)^{A}/((\mathbb{C}^*)^2\times(\mathbb{C}^*))$. It can be easily seen (and contained in [A]) that the quotient torus (which controls the deformations) is the cokernel of the map $((\mathbb{C}^*)^2)^v \rightarrow (\mathbb{C}^*)^e$.

It agrees with the fact that the space of first-order deformations of the directions modulo the subspace of deformations not destroying too rigid graph is the cokernel of the complex $(\mathbb{R}^2)^v \rightarrow (\mathbb{R})^e$ (and these too complexes can be organised into a commutative diagram via exponentiation). So, from the first-order deformation of the directions one obtains the $1$-parametric subgroup in the normal quotient torus.

Being deformed along this (general enough) $1$-parametric subgroup, $L$ tends to the union of boundary components. The boundary component is being tended to by this $1$-parametric subgroup if and only if it is contained in the perturbation.

\note In the complex version these components might have some integer coefficients which are to be calculated.

\section{Questions}

\question What is the complex analogue of the presented structure? From the dimension one we know that it should be a formal series in some variables. Will these variables correspond to the cohomology of $X$ or the space of curves on $X$? How to do it?

\question What should one do if the set $A$ has triples of points on a same line? Even if the perturbation can be done in such a way that avoids too rigid graphs (and it seems it can), it is still unclear what should be used instead of Leibniz rule. The structure becomes definitely something richer than just a cyclic operad.

\question There is a problem with a formalism of virtual fundamental classes for this theory - the normal sheaf to the surface can have second cohomology, hence the deformation theory is not perfectly obstructed. The work [M] solves this problem (in a way we don't fully understand). Can one use these virtual fundamental classes to prove that structures of the type considered in this paper actually act on any projective variety $X$ with an action of $(\mathbb{C}^*)^2$?

\question It is unclear how to do the same formalism of perturbations in the dimension greater than $2$.

\end{document}